\documentclass[reqno]{cpde}

\newtheorem{theorem}{Theorem}%[section]
\newtheorem{lemma}[theorem]{Lemma}
\newtheorem{corollary}[theorem]{Corollary}

\theoremstyle{remark}
\newtheorem{remark}{Remark}%[section]
\newtheorem{step}{Step}

\numberwithin{equation}{section}

\newcommand*  {\e}{{\mathrm e}}
\renewcommand*{\d}{{\mathrm d}}

\newcommand*  {\R} {{\mathbb R}}
\newcommand*  {\N} {{\mathbb N}}

\newcommand*{\norm}[3][{\vphantom 1}]{\lVert #2 \rVert_{#3}^{#1}}

\newcommand*{\Space}[3]{{#1}_{{\vphantom 1}#2}^{{\vphantom +}#3}}
\renewcommand*{\l}[2][]{\Space{L}{#1}{#2}}

\renewcommand*{\c}[2][]{\Space{C}{#1}{#2}}

\allowdisplaybreaks

\begin{document}

\title[The vortex blob method as a second-grade fluid]
{The vortex blob method as a second-grade non-Newtonian fluid}

\author[M. Oliver]{Marcel Oliver}
\address[M. Oliver]{Mathematisches Institut \\
Universit\"at T\"ubingen \\
72076 T\"ubingen \\
Germany}
\email{oliver@member.ams.org}
\author[S. Shkoller]{Steve Shkoller}
\address[S. Shkoller]{Department of Mathematics\\
University of California \\
Davis, CA 95616 \\
USA}
\email{shkoller@math.ucdavis.edu}

\year{June 15, 1999; current version October 14, 1999}

\begin{abstract}
We show that a certain class of vortex blob approximations for ideal
hydrodynamics in two dimensions can be rigorously understood as
solutions to the equations of second-grade non-Newtonian fluids with
zero viscosity, and initial data in the space of Radon measures
${\mathcal M}({\mathbb R}^2)$.  The solutions of this regularized PDE,
also known as the averaged Euler or Euler-$\alpha$ equations, are
geodesics on the volume preserving diffeomorphism group with respect
to a new weak right invariant metric.  We prove global existence of
unique weak solutions (geodesics) for initial vorticity in ${\mathcal
M}({\mathbb R}^2)$ such as point-vortex data, and show that the
associated coadjoint orbit is preserved by the flow.  Moreover,
solutions of this particular vortex blob method converge to solutions
of the Euler equations with bounded initial vorticity, provided
that the initial data is approximated weakly in measure, and the total
variation of the approximation also converges.  In particular, this
includes grid-based approximation schemes of the type that are usually
used for vortex methods.
\end{abstract}

\maketitle

\section{Introduction}

The starting point of our investigation is the somewhat surprising
fact that the equations of motion for an inviscid non-Newtonian fluid
of second grade, and Chorin's vortex blob algorithm with a particular
choice of cut-off or blob function are, at least formally, equivalent.

The velocity field $u=u(x,t)$ of a second grade fluid, under the
assumptions of observer objectivity and material frame-indifference,
satisfies the unique equation
\begin{subequations}
  \label{e.2grade}
\begin{gather}
    (1-\alpha^2 \Delta) \partial_t u +
    u \cdot \nabla(1-\alpha^2\Delta)u -
    \alpha^2 \, (\nabla u)^t \cdot \Delta u
    = -\operatorname{grad} p \,, 
      \label{e.2grade.a} \\
    \operatorname{div} u = 0 \,,
      \label{e.2grade.b} \\
    u(0) = u_0 \,,
      \label{e.2grade.c} 
\end{gather}
\end{subequations}
where $p = p(x,t)$ is the pressure function which is determined
(modulo constants) by the velocity field.  See \cite{DF} and
references therein for a discussion of the constitutive theory of
second grade fluids, and \cite{CO,CG} for well-posedness of the
viscous second-grade fluid equations.  In this context, the constant
$\alpha>0$ is a material parameter which represents the elastic
response of the fluid.

In two dimensions, taking the curl of equation \eqref{e.2grade.a} and
setting $q=(1-\alpha^2\Delta) \operatorname{curl}_{\operatorname{2D}}
u$ yields the vorticity form
\begin{subequations}
  \label{e.alpha-vorticity}
\begin{gather}
  \partial_t q + u \cdot \operatorname{grad} q = 0 \,,
    \label{e.alpha-vorticity.a} \\
  u = K^\alpha * q \,,
    \label{e.alpha-vorticity.b} \\
  q(0)=q_0 \,, 
\end{gather}
\end{subequations}
where $q=q(x,t)$ is called the \emph{potential vorticity}, and
$K^\alpha$ is the integral kernel of the inverse of $(1-\alpha^2
\Delta) \operatorname{curl}_{\operatorname{2D}}$, defined so that the
divergence condition \eqref{e.2grade.b} is satisfied.  

When $\alpha$ is interpreted as a length scale, \eqref{e.2grade} or
\eqref{e.alpha-vorticity} are known as the averaged Euler or
Euler-$\alpha$ equations \cite{HMR} which model the large scale flow
(spatial scales larger than $\alpha$) of an ideal incompressible
fluid.  Their analysis and rich geometry has recently received much
attention \cite{MRS,S1,S2}.  In particular, solutions of
\eqref{e.2grade} on an $n$-dimensional Riemannian manifolds $(M,g)$
arise as geodesic flow on the group of $H^s$-class volume preserving
diffeomorphisms ${\mathcal D}_\mu^s$ provided $s>(n/2)+1$ with respect
to a new weak right invariant metric, given at the identity element $e
\in {\mathcal D}_\mu^s$ by
\begin{equation}\label{metric}
\langle u, v \rangle_e = (u,v)_{\l2} + 
         2 \alpha^2 (\operatorname{Def}u, \operatorname{Def}v)_{\l2} 
\end{equation}
where $\operatorname{Def}u = (\nabla u + \nabla u^t)/2$.  Thus,
following the program of Arnold \cite{A} and Ebin-Marsden \cite{EM},
local-in-time well-posedness of classical solutions is a direct
consequence of the existence of $C^\infty$ geodesics of $\langle
\,\cdot\,, \,\cdot\, \rangle$ on ${\mathcal D}_\mu^s$.

The vortex blob method was introduced by Chorin \cite{C} as a
regularization of the point vortex algorithm for ideal hydrodynamics,
and can be understood as follows.  Consider the vorticity form of the
Euler equations on $\R^2$,
\begin{subequations}
  \label{e.euler}
\begin{gather}
  \partial_t \omega + u \cdot \operatorname{grad} \omega = 0 \,,
    \label{e.euler.a} \\
     u = K * \omega, \label{e.euler.b} \\
  \omega(0) = \omega_0 \,.
\end{gather}
\end{subequations}
Here $K(x,y) = 1/(2 \pi) \, \nabla^\perp \log \left| x-y \right|$ and
$\omega = \omega (x,t)$ is the physical vorticity of the flow. When
the velocity field is sufficiently regular---$u$ is at least
continuous in $t$ and quasi-Lipschitz in $x$, uniformly over finite
intervals of time---we may define the Lagrangian flow map
$\eta_t=\eta( \,\cdot\,, t)$ by
\begin{equation}
  \label{e.flow-map}
  \partial_t \eta(x,t) = u(\eta(x,t), t) \,,
\end{equation}
or equivalently by
\begin{equation}
  \partial_t \eta_t = u_t \circ \eta_t \,.
\end{equation}
For each $t$, the map $\eta( \,\cdot\,, t)$ is in ${\mathcal G}$, the
group of all homeomorphisms $\phi$ of ${\mathbb R}^2$ which preserve
the Lebesgue measure. The pointwise conservation of vorticity under
the Euler flow is thus expressed by $\omega_t \circ \eta_t =
\omega_0$; combining \eqref{e.flow-map}, \eqref{e.euler.b}, and the
initial condition $\eta( \,\cdot\,, 0) = e$, we obtain the ODE
\begin{equation}
  \label{e.lagrangian}
  \partial_t \eta (x,t)
  = \int_{\R^2} K(\eta(x,t),\eta(y,t)) \, \omega_0(y) \, \d y \,. 
\end{equation}
Letting $\delta$ denote the Dirac measure and substituting the point
vortex ansatz
\begin{equation}
  \label{e.ansatz}
  \omega (x,t) = \sum_{n=1}^N \Gamma_i \, \delta (x - x_i(t)) \,,
\end{equation}
into \eqref{e.lagrangian}, we obtain a finite dimensional system of
ordinary differential equations for the vortex centers $x_1, \dots,
x_N$.  However, the induced velocity field has $1/|x|$-type
singularities at the vortex centers.  Hence, the point vortex system
is neither numerically well-behaved (the exact solution of the point
vortex system may even collapse in a finite time for small sets of
initial data \cite{MP94}), nor does it approximate physically relevant
velocity fields very well.

The idea of the vortex blob method is to smooth the Dirac measure by a
\emph{cut-off} or \emph{blob function} $\chi$ that decays at infinity
and whose mass is mostly supported in a disc of diameter $\alpha$.
This leads to the following equation for the Lagrangian flow:
\begin{subequations}
  \label{e.alpha-lagrangian}
\begin{gather}
  \partial_t \eta^\alpha (x,t)
  = \int_{\R^2} K^\alpha(\eta^\alpha (x,t), \eta^\alpha (y,t)) \,
    \omega_0(y) \, \d y \,,
  \label{e.alpha-lagrangian.a}   
\intertext{where}  
  K^\alpha = \nabla^\perp G^\alpha \,, \\
  -\Delta G^\alpha(x,y) = \chi^\alpha(|x-y|)
  \equiv \frac1{\alpha^2} \chi \biggl( \frac{|x-y|}\alpha \biggr) \,.
\end{gather}
\end{subequations}
Many researchers have investigated the convergence properties of this
scheme \cite{BM82,BM85,C88,H87,MP82}.  In particular, for certain
smooth cut-off functions, such as Bessel functions, the order of
accuracy with respect to the regularization parameter $\alpha$ depends
only on the smoothness of the Euler flow (``infinite order
accuracy'').

It is now easy to see that the equation of a second grade fluid
\eqref{e.alpha-vorticity} and the vortex blob method coincide when
$\chi(x) = K_0(x)$, where $K_0$ is the modified Bessel function of the
second kind which is the Green's kernel for the operator $(1-\Delta)$.
Thus far, this relationship has only been formally established, as it
remains to be proven that the point-vortex ansatz \eqref{e.ansatz}
makes sense as data for the PDE \eqref{e.alpha-vorticity}; moreover,
it is not \emph{a priori} clear if solutions to the vortex blob method
converge to true Euler solutions (in the sense of PDE).  Our results
are the following.

We show that the Lagrangian flow formulation of the blob method
\eqref{e.alpha-lagrangian} with $K_0$ cut-off function is well-posed
for initial potential vorticities $q_0$ in $\mathcal{M}(\R^2)$, the
space of Radon measures on $\R^2$.  In particular, this includes
point-vortex initial data.  Such a result does not hold for the Euler
equations, where the flow map of the point vortex system
\eqref{e.lagrangian} is not known to be well-defined.

This result allows us to rigorously classify the co-adjoint orbits
characterized by point-vortex initial data.  Let us explain what we
mean by this.  The configuration space for ideal incompressible
hydrodynamics is the volume-preserving diffeomorphism group, and for
$s>(n/2)+1$, ${\mathcal D}_\mu^s$ is a $C^\infty$ Hilbert manifold and
a smooth topological group.  While ${\mathcal D}_\mu^s$ is not a Lie
group (left composition and inversion are only $C^0$ and the group
exponential map does not cover a neighborhood of the identity), it
behaves similar to a Lie group, because of the smooth properties of
the Riemannian exponential map (see \cite{EM} and \cite{S1,S2}).  The
Eulerian phase space for the fluid motion is the single fiber
$T_e{\mathcal D}_\mu^s$ consisting of $H^s$-class divergence free
vector fields on the fluid container, and this vector space can be
formally thought of as the ``Lie algebra'' of ${\mathcal D}_\mu^s$.
The cotangent space at the identity is given by $H^s(\Lambda^1)/d
H^s(\Lambda^0)$, the $H^s$-class differential $1$-forms modulo exact
$0$-forms.  Using the fact that the exterior derivative $d\colon
T^*_e{\mathcal D}_\mu^s \rightarrow H^{s-1}(\Lambda^2)$ is an
isomorphism, and the fact that we may identify $H^{s-1}(\Lambda^2)$
with $H^{s-1}(\Lambda^0)$, the role of the dual of the ``Lie algebra''
for 2D hydrodynamics is played by the $H^{s-1}$-class vorticity
functions.  The representation of ${\mathcal D}_\mu^s$ on this ``Lie
algebra'' is provided by the co-Adjoint action, so that for $\eta\in
{\mathcal D}_\mu^s$ and $\omega \in H^{s-1}(\Lambda^0)$,
$\operatorname{Ad}^*_\eta(\omega) = \omega \circ \eta$, and the
invariance of the co-Adjoint orbit is merely the pointwise
conservation of vorticity which is fundamental to 2D hydrodynamics.
If one temporarily ignores the topology and works formally, then it is
possible to classify certain interesting and important co-Adjoint
orbits.  Specifically, it is a result of Marsden and Weinstein
\cite{MW} that point-vortex initial data \eqref{e.ansatz} define the
co-Adjoint orbit on which point-vortex dynamics evolve.  This is
clearly a formal result as Dirac measures are not elements of
$H^{s-1}$ for $s>2$; consequently, the problem is to supply a
candidate topology for the ``Lie algebra'' which is general enough to
contain the Dirac measures, and weaken the regularity of the
configuration space so that its ``representation'' is well-defined.
In doing so, one can establish a rigorous classification of the orbit.
By using ${\mathcal G}$ for the configuration space and ${\mathcal
M}(\R^2)$ for the ``Lie algebra,'' and by defining a new notion of
{\it weak} co-adjoint action which coincides with the notion of a weak
solution, we are able to establish the orbit classification for
point-vortex initial data, and prove that our particular vortex blob
method leaves such weak co-adjoint orbits invariant.

Finally, we consider the matter of greatest practical importance: the
convergence of solutions of the vortex blob method to solutions of the
Euler equations as the blob diameter $\alpha \rightarrow 0$.  We prove
this convergence result under the rather mild assumption that the
initial Euler vorticity field $\omega_0$ is continuous with compact
support and is approximated on its support by a sequence of weakly
converging measures in $\mathcal{M}(\R^2)$ that have uniformly bounded
total variation. (The restriction to compact support will be replaced
by weaker assumptions on the decay at infinity.)

The precise statements of our results are as follows.

\begin{theorem} \label{t.well-posedness}
  For initial data $q_0\in {\mathcal M}(\R^2)$, there exists a unique
  global weak solution to \eqref{e.alpha-vorticity} with
\begin{equation}
 \eta^\alpha \in \c1 (\R; {\mathcal G}) \,, \quad
 u^\alpha\in \c0(\R; C^0_{\mathrm{div}}(\R^2, \R^2)) \,, 
 \quad \text{and} \quad 
 q \in C^0(\R; {\mathcal M}(\R^2)) \,,
\end{equation}
where the subscript ${\mathrm {div}}$ denotes divergence-free.  As a
consequence, the co-Adjoint action $\operatorname{Ad}^*_{\eta}(q)$ and
the weak co-adjoint action $\operatorname{w-ad}^*_{u}(q)$ are
conserved.
\end{theorem}

\begin{remark}
  The solution that we construct may not necessarily have finite
  energy, i.e., the velocity field $u^\alpha$ may not be in $\l2$.
  None of our results, however, relies on energy type estimates.
  Furthermore, as is the case for the Euler equations, the initial
  potential vorticity can be decomposed into a radially symmetric and
  a mean-zero part, with a corresponding velocity field $u$ in the
  affine space $u_{\operatorname{stationary}} + \l2 (\R^2,\R^2)$.  For
  details see DiPerna and Majda \cite{DM87}.
\end{remark}

\begin{remark}
  An immediate consequence of the uniqueness of the solution and the
  time-reversibility of the equation is that the vortex blob system
  cannot collapse in finite time, i.e., two or more vortex centers
  cannot merge into one in finite time.  For non-regularized Euler
  vortex dynamics, on the other hand, it is known that vortex collaps
  occurs on small sets of initial configurations \cite{MP94}.
\end{remark}

\begin{remark}
  The kernel $K^\alpha$ which corresponds to a second grade fluid is
  the least regular kernel (modulo possible sub-logarithmic
  corrections) for which uniqueness of point vortex solutions can be
  shown.  An equivalent uniqueness result based on Sobolev space
  methods, and for bounded domains is given in \cite{HKOT}.
\end{remark}

\begin{theorem} \label{t.convergence}
Let $\eta$ be the flow map of the Euler equation \eqref{e.lagrangian}
with initial vorticity $\omega_0 \in \l1(\R^2) \cap L^\infty(\R^2)$. Suppose
that $\omega_0$ is approximated by a sequence of measures $q_0^n$ in
${\mathcal M}(\R^2)$ such that $q_0^n \rightharpoonup \omega_0$ weakly
in $\mathcal{M}(\R^2)$ and $\| q_0^n \|_{\mathcal M} \to \| \omega_0
\|_{\l1}$.  Then for every $T>0$, there exists a sequence
$\{\alpha_n\}$ converging to zero as $n \to \infty$ such that when
$\eta^{\alpha_n}$ denotes the flow map of the vortex method with
$\alpha = \alpha_n$ and initial data $q_0^n$,
\begin{equation}
  \lim_{n \to \infty}
  \sup_{t \in [0,T]} \sup_{x \in \R^2}
    \bigl| \eta^{\alpha_n} (x,t) - \eta (x,t) \bigr| = 0 \,.
\end{equation}
\end{theorem}

\begin{remark}
  The idea of analyzing the vortex method as a PDE posed on some space
  of distributions was already used by Marchioro-Pulvirenti \cite{MP82}
  and Cottet \cite{C88}.  Cottet's result
  requires stronger assumptions on the cut-off function and
  hence smoothing kernel, stronger regularity assumptions on the
  underlying Euler flow, and his approximation of the Euler initial
  vorticity field required a uniform grid.  The trade-off, however, is
  that these more stringent constraints give an improved (algebraic)
  convergence in $\alpha$.  There is, in general, a trade-off between
  the order of convergence on the one hand, and the assumptions placed
  on $K^\alpha$, $\omega_0$, and the approximation at time $t=0$ on
  the other.  The result which, to our knowledge, comes closest to
  Theorem~\ref{t.convergence} is given in Marchioro and Pulvirenti
  \cite{MP82}.  The authors, however, assume that $K^\alpha$ is
  Lipschitz, which again excludes kernels corresponding to the
  equations of second grade fluids.
\end{remark}

\begin{remark}
Stronger results can be proved for kernels $K^\alpha$ with a higher
degree of smoothing.  For example, by replacing $1-\alpha^2 \Delta$
with $(1-\alpha^2\Delta)^s$, one obtains a hierarchy of
regularizations of the Euler equations which coincide with geodesic
flow on the volume-preserving diffeomorphism group with respect to the
$H^s$ metric.  Other choices of $K^\alpha$ may introduce non-local
pseudo-differential operators into equation \eqref{e.alpha-vorticity},
but the analysis can still proceed as before.
\end{remark}

\begin{remark}
In three dimensions, the formal connection between second grade fluids
and particular vortex filament methods still holds and is the subject
of a forthcoming article.  In this setting, one looks at the set of
vorticity distributions of the following form. Let $\gamma$ be a curve
in $ {\R}^3 $ extending to infinity in both directions, and let $
\delta_\gamma $ be the Dirac distribution given by integration along
$\gamma$ with respect to arc length.  Let $ \omega_\gamma $ be the
2-form along $\gamma$ defined by $ i_T dx \wedge dz$, where $T$ is the
unit tangent vector to $\gamma$.  Then if $\Gamma$ is any constant, $
\Gamma \omega _\gamma \delta_\gamma $ is the vorticity corresponding
to $\gamma$ with strength $ \Gamma $. See \cite{MW}.
\end{remark}

\begin{remark}
As we described above, for $s>(n/2)+1$, local well-posedness follows
form the existence of unique $\c\infty$ geodesics $\dot\eta(t)$ on
${\mathcal D}_\mu^s$ with respect to the right invariant metric
$\langle \cdot, \cdot \rangle$ defined in (\ref{metric}), with initial
conditions $\eta(0)=e$ and $\dot\eta(0)=u_0$. In working with the
geodesic flow ${\dot{\eta}}(t)$, one obtains $C^\infty$ evolution in
the tangent bundle $T{\mathcal D}_\mu^s$ and $C^\infty$ dependence on
initial data, while the projected evolution curve $u(t) =
{\dot{\eta}}(t) \circ \eta(t)^{-1}$ in the single fiber of the tangent
bundle $T_e {\mathcal D}_\mu^s$ -- which plays the role of the
Eulerian phase space -- has only $C^0$ smoothness, and $C^0$
dependence on the initial velocity field.  In the case that the
manifold has a smooth boundary, there are three new subgroups of
${\mathcal D}_\mu^s$ which are in one-to-one correspondence with the
classical Dirichlet, Neumann, and mixed elliptic boundary value
problems in the sense that elements of the ``Lie algebras'' of these
three subgroups satisfy those boundary conditions.  Hence, geodesic
flow of $\langle \cdot, \cdot \rangle$ on these three subgroups gives
the solutions of (\ref{e.2grade}) with no-slip, free-slip, and mixed
boundary conditions.
\end{remark}

\section{Kernel estimates}

The crucial ingredients for the proof of our theorems are
quasi-Lipschitz estimates on the Euler kernel $K$ and the regularized
kernel $K^\alpha$.  For $x \in \R^2$ we define the function
\begin{equation}
  \label{e.phi}
  \varphi (x) =
  \begin{cases}
    |x| \, (1- \ln {|x|})  & \text{for } |x| < 1 \,, \\
    1                      & \text{for } |x| \geq 1 \,.
  \end{cases}
\end{equation}
\begin{lemma} \label{l.quasi-lipshitz.1}
  For $\omega \in \l1(\R^2) \cap \l\infty(\R^2)$,
\begin{equation}
    \int_{\R^2} \bigl| K(x,y) - K(x',y) \bigr| \, |\omega(y)| \, \d y 
      \leq c \, \varphi (x-x') \,
        (\norm{\omega}{\l1} + \norm{\omega}{\l\infty}) \,.
\end{equation}
\end{lemma}      

The proof is standard and can be found, for example, in McGrath
\cite{M68}.  Somewhat less standard is the following estimate, still
for the Euler kernel, which is similar to estimates in Benedetto
\emph{et al.}\ \cite{BMP93}.

\begin{lemma} \label{l.quasi-lipshitz.2}
Let $\omega \in \l1(\R^2) \cap \l\infty(\R^2)$ and let $\phi$ be an
area preserving measurable transformation on $\R^2$.  Then
\begin{equation}
  \biggl|
    \int_{\R^2}
      \bigl( K(x,y) - K(x,\phi(y)) \bigr) \, \omega(y) \, \d y
  \biggr|
  \leq c \, \sup_{x \in \R^2} \varphi \bigl( x - \phi(x) \bigr) \,
       \bigl(
         \norm{\omega}{\l1} + \norm{\omega}{\l\infty}
       \bigr) \,.
  \label{e.k-estimate.2}
\end{equation}
\end{lemma}

\begin{proof}
Set $r = \sup_{x} \lvert x - \phi(x) \rvert$; as in the proof of
Lemma~\ref{l.quasi-lipshitz.1}, the interesting case is when $r<1$.
We split the integral in \eqref{e.k-estimate.2} into two parts.
First, consider
\begin{align}
  \int_{\lvert x-y \rvert \leq 2r}
  & \bigl| K(x,y) - K(x,\phi(y)) \bigr| \, |\omega(y)| \, \d y
           \notag \\
  & \leq \frac1{2\pi} \int_{\lvert x-y \rvert \leq 2r}
           \frac{|\omega(y)|}{|x-y|} \, \d y +
         \int_{\lvert x-y \rvert \leq 2r}
           \frac{|\omega(y)|}{|x-\phi(y)|} \, \d y \notag \\
  & \leq \frac1{2\pi} \int_{\lvert x-y \rvert \leq 2r}
           \biggl(
             \frac1{|x-y|} + \frac1{|x-\phi(y)|}
           \biggr) \d y \,
           \norm{\omega}{\l\infty} \notag \\
  & \leq \frac1{\pi} \int_{\lvert x-y \rvert \leq 2r}
             \frac1{|x-y|}\, \d y \,
           \norm{\omega}{\l\infty}
    \equiv c \, r \, \norm{\omega}{\l\infty} \,.
\end{align}
The last inequality above holds because $\phi$ is area preserving, and
among all such transformations, the symmetric map
$\phi=e$ maximizes the integral over $|x-\phi(y)|^{-1}$.

Next, we consider the case when $\lvert x-y \rvert \geq 2r$.  Observe
that
\begin{equation}
       \lvert x - \phi(y) \rvert
  \geq \lvert x - y \rvert - \lvert y - \phi(y) \rvert
  \geq \lvert x - y \rvert - r
  \geq \tfrac12 \, \lvert x - y \rvert \,,
\end{equation}
so that
\begin{align}
  \int_{\lvert x-y \rvert \geq 2r}
  & \bigl| K(x,y) - K(x,\phi(y)) \bigr| \, |\omega(y)| \, \d y
           \notag \\
  & \leq \frac1{2\pi} \int_{\lvert x-y \rvert \geq 2r}
           \frac{| y - \phi(y)|}{|x-y| \, |x - \phi(y)|} \,
           |\omega(y)| \, \d y \notag \\
  & \leq \frac1{\pi} \int_{\lvert x-y \rvert \geq 2r}
           \frac{r}{|x-y|^2} \,
           |\omega(y)| \, \d y \notag \\  
  & \leq \frac{r}\pi
           \left(
             \int_{2r \leq \lvert x-y \rvert \leq 2}
               \frac{|\omega(y)|}{|x-y|^2} \, \d y
             + \int_{\lvert x-y \rvert \geq 2}
               \frac{|\omega(y)|}{|x-y|^2} \, \d y
           \right) \notag \\
  & \leq \frac{r}\pi
           \left(
             \int_{2r}^2 \frac{\d \rho}\rho \,
               \norm{\omega}{\l\infty} +
             \frac14 \, \norm{\omega}{\l1}
           \right) \notag \\
  & \leq c \, \varphi(r) \,
           \bigl(
             \norm{\omega}{\l1} + \norm{\omega}{\l\infty}
           \bigr) \,.
\end{align}
By combining the two estimates we complete the proof.
\end{proof}

Finally, we give the corresponding result for the vortex method kernel.

\begin{lemma} \label{l.quasi-lipshitz.3}
  There exists a constant $c_2$ which is independent of $\alpha$, such
  that
\begin{equation}
  \sup_{y \in \R^2}
       \bigl| K^\alpha (x,y) - K^\alpha(x',y) \bigr|
  \leq \frac{c_2}\alpha \, 
       \varphi \Bigl( \frac{x-x'}\alpha \Bigr) \,.
\end{equation}
\end{lemma}

\begin{proof}
  Note that on $\R^2$, $K^\alpha(x,y) = K^\alpha(|x-y|) = \nabla^\bot
  G^\alpha (|x-y|)$, where
\begin{equation}
  G^\alpha(r) = \frac1{2\pi} \, K_0 \Bigl( \frac{r}\alpha \Bigr)
                + \frac1{2\pi} \, \ln r 
\end{equation}  
and $K_0$ denotes the zero order modified Bessel function of the
second kind \cite{AS64}.  For simplicity, we take $\alpha=1$ and
compute
\begin{align}
  \frac{\d G^\alpha}{\d r} (r)
    & = \frac1{2\pi} \left( \frac1r - K_1 (r) \right)
      = \frac1{4\pi} \, r \, \ln r + O (r) \,,  \\
  \frac{\d^2 G^\alpha}{\d r^2} (r)
    & = \frac1{2\pi} \left( \frac1{r^2} + K_0 (r)
        + \frac1r \, K_1 (r) \right)
      = - \frac1{4 \pi} \, \ln r + O(1) \,,
    \label{e.grr-expansion}
\end{align}
as $r \to 0$.  Set $r \equiv |x-x'|$ and assume, without loss of
generality as $K^\alpha$ is bounded, that $r<1$.

If $|x-y| < 2r$, then $|x'-y| \leq |x'-x| + |x-y| < 3r$, so that
\begin{align}
  \bigl| K^\alpha (x,y) - K^\alpha (x',y) \bigr|
  & \leq \bigl| \nabla^\bot G^\alpha (|x-y|) \bigr| +
         \bigl| \nabla^\bot G^\alpha (|x-y|) \bigr|     \notag \\
  & \leq \biggl| \frac{\d G^\alpha}{\d r} (|x-y|) \biggr|+
         \biggl| \frac{\d G^\alpha}{\d r} (|x-y|) \biggr|   \notag \\
  & \leq \frac2\pi \, r \, \lvert \ln r \rvert + O(r) \,.
\end{align}
Since $\d G^\alpha/{\d r}$ is continuous and decays at infinity, this
implies a bound of the form
\begin{equation}
  \label{e.2}
  \bigl| K^\alpha (x,y) - K^\alpha(x',y) \bigr|
    \leq c \, \varphi(|x-x'|) \,.
\end{equation}

If, on the other hand, $|x-y| \geq 2r$, we use the mean value theorem
to estimate
\begin{align}
  \bigl| K^\alpha (x,y) - K^\alpha (x',y) \bigr|
  & \leq \sup_{x'' \in B(x,r)} |\nabla K^\alpha (x'',y)| \, |x-x'|
         \notag \\
  & \leq \biggl| \frac{\d^2 G^\alpha}{\d r^2} (|x''-y|) \biggr| \, r 
         \notag \\
  & \leq \frac1{4\pi} \, r \, \ln r + O(r) \,,
\end{align}
which again implies a bound of the form \eqref{e.2}. In the last step
we have used \eqref{e.grr-expansion} in conjunction with $|x''-y| >
r$.

To recover the scaling of the estimate in $\alpha$, divide
\eqref{e.2} by $\alpha$, rescale $x$, $x'$, and $y$ by $\alpha^{-1}$,
and note that $K^\alpha(r) = K^{\alpha=1}(r/\alpha)/\alpha$.
\end{proof}

\begin{corollary} 
  For $q \in \mathcal{M}(\R^2)$,
\begin{equation}
    \int_{\R^2}
        \bigl| K^\alpha (x,y) - K^\alpha (x',y) \bigr| \, 
        |q(y)| \,  \d y 
      \leq \frac{c_2}\alpha \,
        \varphi \Bigl( \frac{x-x'}\alpha \Bigr) \,
        \norm{q}{\mathcal{M}}  \,.
\end{equation}
\end{corollary}

\section{Well-posedness}

We can now prove the existence of unique, global, weak solutions to
the Lagrangian flow equation \eqref{e.alpha-lagrangian}.

\begin{proof}[Proof of Theorem~\ref{t.well-posedness}]
Due to the quasi-Lipschitz condition for $K^\alpha$, we can adopt the
method that Kato developed for the Euler equations in \cite{K67}, by
simply replacing the kernel estimates in $\l1$ by the corresponding
estimates in $\l\infty$.  Our presentation follows to some extent that
of Marchioro and Pulvirenti \cite{MP94}.

For simplicity, we assume $\alpha = 1$ throughout this proof. We
introduce a sequence of approximate solutions
\begin{subequations}
  \label{e.approx}
  \begin{gather}
    \partial_t \eta^n(x,t) = u^n (\eta^n(x,t),t) \,,
                                                \label{e.approx.a} \\   
    \eta^n(x,0) = x \,,                         \label{e.approx.b} \\
    \eta^0(x,t) = x \,,                         \label{e.approx.c} \\
    q^n (\eta^n(x,t),t) = q_0(x) \,,            \label{e.approx.d} \\
    u^n(x,t) = \int_{\R^2}
      K^\alpha(x,y) \, q^{n-1}(y,t) \, \d y \,, \label{e.approx.e}
  \end{gather}
\end{subequations}
for $n \in \N$. The proof now proceeds in several steps.

\begin{step}
Prove that $\eta^n \in \c1 ((0, \infty); {\mathcal G})$ for every $n
\in \N$.
\end{step}

We proceed inductively.  Notice that for every $n$ the vector field
$u^n$ is quasi-Lipschitz in space and continuous in time.  This is a
consequence of Lemma~\ref{l.quasi-lipshitz.3} as
\begin{align}
  \bigl| u^n (x,t) - u^n(x',t) \bigr|
  & =    \biggl|
           \int_{\R^2}
             \Bigl[
               K^\alpha (x, \eta^{n-1}(y,t)) -
               K^\alpha (x', \eta^{n-1}(y,t))
             \Bigr] \,
             q_0(y) \, \d y
         \biggr| \notag \\
  & \leq c\, \varphi (x-x') \, \norm{q_0}{\mathcal{M}} \,,
\intertext{and}
\bigl| u^n (x,t) - u^n(x,t') \bigr|
  & =    \biggl|
           \int_{\R^2}
             \Bigl[
               K^\alpha (x, \eta^{n-1}(y,t)) -
               K^\alpha (x, \eta^{n-1}(y,t'))
             \Bigr] \,
             q_0(y) \, \d y
         \biggr| \notag \\
  & \leq \sup_{y \in \R^2}
           \bigl|
               K^\alpha (x, \eta^{n-1}(y,t)) -
               K^\alpha (x, \eta^{n-1}(y,t'))
           \bigr| \, 
           \norm{q_0}{\mathcal{M}} \notag \\
  & \leq c \, \sup_{y \in \R^2} \varphi
           \bigl(
               \eta^{n-1}(y,t) - \eta^{n-1}(y,t')
           \bigr) \, 
           \norm{q_0}{\mathcal{M}} \notag \\
  & \leq c \, \sup_{y \in \R^2} \sup_{x \in [t,t']} \varphi
           \bigl(
              | \dot \eta^{n-1}(y,s) | \, |t-t'|
           \bigr) \, 
           \norm{q_0}{\mathcal{M}} \notag \\         
  & =    c \, \sup_{y \in \R^2} \sup_{x \in [t,t']} \varphi
           \bigl(
              | u^{n-1}(y,s) | \, |t-t'|
           \bigr) \, 
           \norm{q_0}{\mathcal{M}} \,.
\end{align}
This implies uniform continuity in time, because $u^n$ is bounded for
every $n$:
\begin{equation}
  |u^n(x,t)| 
  =      \biggl|
           \int_{\R^2} K^\alpha (x, \eta^{n-1}(y,t)) \,
             q_0(y) \, \d y
         \biggr|
  \leq   \sup_{y \in \R^2} 
           \bigl| K(x,y) \bigr| \, \norm{q_0}{\mathcal{M}}
  \equiv c \, \norm{q_0}{\mathcal{M}} \,.
  \label{e.u-bound}
\end{equation}
Since $u^n$ is continuous in time and quasi-Lipschitz in space, the
vector field generates a local flow $\eta^n \in \c1([0,T); C(\R^2))$
for some $T>0$---see, e.g., Chapter~2, Lemma~3.2 in Marchioro and
Pulvirenti~\cite{MP94}.  Because of the global bound
\eqref{e.u-bound}, the right side of \eqref{e.approx.a} is bounded and
the flow exists globally in time.

\begin{step} \label{s.2}
Show that there exists a limiting flow map $\eta \in C ((0, \infty);
{\mathcal G})$.
\end{step}

We first prove that the sequence $\eta^n$ is Cauchy in $C([0, T];
{\mathcal G})$ for some $T>0$.  To simplify notation, we shall drop
the explicit time dependence of $u$ and $\eta$, and estimate
\begin{align}
  \bigl| \eta^n & (x,t) - \eta^{n-1}(x,t) \bigr| \notag \\
  & \leq \int_0^t
         \bigl| u^n (\eta^n) - u^{n-1} (\eta^{n-1}) \bigr| \d s
         \notag \\
  & \leq \int_0^t
         \biggl| 
           \int_{\R^2}
             \Bigl[
               K^\alpha (\eta^n(x),\eta^{n-1}(y)) -
               K^\alpha (\eta^{n-1}(x),\eta^{n-1}(y))
             \Bigr]
             q_0(y) \, \d y
         \biggr| \d s \notag \\
  & \quad + \int_0^t
         \biggl| 
           \int_{\R^2}
             \Bigl[
               K^\alpha (\eta^{n-1}(x),\eta^{n-1}(y)) -
               K^\alpha (\eta^{n-1}(x),\eta^{n-2}(y))
             \Bigr]
             q_0(y) \, \d y
         \biggr| \d s \notag \\
  & \leq c \int_0^t
         \varphi \bigl(
                   \eta^n(x) - \eta^{n-1}(x)
                 \bigr) \,
         \d s \, \norm{q_0}{\mathcal{M}} \notag \\
  & \quad + c \int_0^t
         \varphi \bigl(
                   \eta^{n-1}(x) - \eta^{n-2}(x)
                 \bigr) \,
         \d s \, \norm{q_0}{\mathcal{M}} \,.
\end{align}    
By taking the supremum over $x$ on both sides, we obtain
\begin{multline}
  \sup_{x \in \R^2}
    \bigl| \eta^n (x,t) - \eta^{n-1}(x,t) \bigr| \\ 
  \leq c \, \norm{q_0}{\mathcal{M}}
       \int_0^t
         \biggl[
           \varphi \biggl(
                     \sup_{x \in \R^2}
                       \bigl| \eta^n(x) - \eta^{n-1}(x) \bigr|
                   \biggr) +
           \varphi \biggl(
                     \sup_{x \in \R^2}
                       \bigl| \eta^{n-1}(x) - \eta^{n-2}(x) \bigr|
                   \biggr)
         \biggr] \d s \,.
\end{multline}
Defining
\begin{equation}
  \label{e.rho}
  \rho^N (t) \equiv
    \sup_{n \geq N} \sup_{x \in \R^2}
      \bigl| \eta^n (x,t) - \eta^{n-1}(x,t) \bigr| \,,
\end{equation}
we can simplify the previous estimate, and obtain
\begin{equation}
  \rho^N (t) \leq c \int_0^t
    \varphi \bigl( \rho^{N-1}(s) \bigr) \d s \,.
\end{equation}
It is well known that this implies
\begin{equation}
  \label{e.rho-limit}
  \lim_{N \to \infty} \rho^N (t) \to 0 \,,
\end{equation}
uniformly on $[0,T]$ for $T$ sufficiently small.  Since $T$ depends
only on $\alpha$ and the $\mathcal{M}$-norm of $q_0$, this result can
be extended to arbitrarily large times.  Thus, the contraction mapping
theorem implies the assertion of Step~\ref{s.2}.

\begin{step}
Show that the Lagrangian flow equation (\ref{e.alpha-lagrangian}) is
satisfied in the limit, and that $\eta \in \c1 (\R; \mathcal{G})$.
\end{step}

We define the limiting potential vorticity $q$ and the limiting
velocity $u$ in the obvious way, and check by direct estimation that
\begin{equation}
  q^n \rightharpoonup q \equiv q_0 \circ \eta^{-1}
  \label{e.q-def}
\end{equation}
weakly in $\mathcal{M}(\R^2)$, and
\begin{equation}
  u^n \to u \equiv K^\alpha * q
\end{equation}
in $C(\R^2)$; both limits are uniform over finite intervals of time.

To prove that $\eta$, $u$, and $q$ solve the limit problem
(\ref{e.alpha-lagrangian}), we consider its integrated version
\begin{align}
  \eta(x,t) & - \int_0^t u(\eta(x,s),s) \, \d s \notag \\
  & = \eta(x,t) - \int_0^t u(\eta(x,s),s) \, \d s
      - \eta^n(x,t) + \int_0^t u^n(\eta^n(x,s),s) \, \d s \notag \\
  & \leq \bigl| \eta(x,t) - \eta^n(x,t) \bigr|
      + \int_0^t
        \bigl| u^n(\eta^n(x,s),s) - u(\eta^n(x,s),s) \bigr| \, \d s
      \notag \\
  & \quad + \int_0^t
        \bigl| u(\eta^n(x,s),s) - u(\eta(x,s),s) \bigr| \, \d s
    \notag \\
  & \to 0 \text{ uniformly in $x$ as } n \to \infty \,.
\end{align}
Thus, the left side must be zero.  Since $u(\eta(x,s),s)$ is
continuous in $x$, we can differentiate with respect to $t$, and find
that $\eta$ satisfies \eqref{e.alpha-lagrangian} and that $\dot\eta$
is in fact continuous.  Due to the time-reversibility of the equation,
the result extends to negative times as well.

Moreover, one can show---first by formal calculation for smooth
function, and then extending by the usual density argument---that the
weak solution $q$ defined through \eqref{e.q-def} satisfies
\begin{equation}
  \int_\R \int_{\R^2}
    \bigl(
       \partial_t \phi + 
       u \cdot \operatorname{grad} \phi
    \bigr) \,
    q \, \d x \, \d t = 0
\end{equation}
for every $\phi \in \c[0]\infty(\R \times \R^2)$.  This shows that
solutions of the vortex method, and hence the equations of
second-grade non-Newtonian fluids, preserve the (weak) co-adjoint
action.

\begin{step}
Prove that the solution is unique.
\end{step}

Uniqueness is shown by a direct estimate on the difference of two flow
maps.  This leads to another log-Gronwall inequality, which can be
treated in the same way as the previous ones; we omit all details.
\end{proof}

\begin{remark}
  The homeomorphisms that we consider have the vector space $\R^2$ as
  the range; we may thus subtract two elements of this class.  For
  homeomorphisms of a compact domain $\Omega$ of $\R^2$, one can
  isometrically embed the set of measure-preserving homeomorphisms of
  $\Omega$ into the vector space $\l2(\Omega, \R^2)$, and take
  differences in this large space.  Similarly, the difference
  $u^\alpha \circ \eta^\alpha - u \circ \eta$ is not an intrinsic
  operation, but rather relies on the trivial identification of vector
  spaces induced by the trivial geometry of $\R^2$.  On the other
  hand, when the configuration space is ${\mathcal D}_\mu^s(M)$, $s>2$
  and $M$ is a compact Riemannian manifold, the map $u^\alpha \circ
  \eta^\alpha$ is an element of the fiber $T_{\eta^\alpha}{\mathcal
  D}_\mu^s$ while $u \circ \eta$ is in $T_\eta{\mathcal D}_\mu^s$;
  thus, in order to compare the two maps, we must parallel transport
  $u\circ \eta$ into $T_{\eta^\alpha}{\mathcal D}_\mu^s$ along the
  Riemannian connection.
\end{remark}

\section{Weak co-adjoint action and reduction}

As we described, classical solutions of the two-dimensional averaged
Euler equations are geodesics on the Hilbert-class volume-preserving
diffeomorphism group ${\mathcal D}_\mu^s$, $s>2$.  We identify the
space of classical vorticity solutions $H^{s-1}(M)$ with the reduced
space $T_e{\mathcal D}_\mu^s = T{\mathcal D}_\mu^s/ {\mathcal
D}_\mu^s$ (symmetry reduction by the massive particle relabeling
symmetry group ${\mathcal D}_\mu^s$ of hydrodynamics), and note that
this space is the union of the ${\mathcal D}_\mu^s$-co-adjoint orbits.

%In the case that $M={\mathbb R}^2$, it is convenient to use the
%$C^{k+\alpha}$ H\"{o}lder-class volume-preserving diffeomorphisms,
%$k\ge 1, \alpha>0$ which we shall denote by ${\mathcal
%H}_\mu^{k+\alpha}$, because this space does not place restrictions of
%decay at $\infty$.

In the case that $M={\mathbb R}^2$, and for the purpose of studying
weak solutions to \eqref{e.alpha-vorticity} we shall substantially
relax the regularity requirements on the configuration space, and use
${\mathcal G}$ in place of ${\mathcal D}_\mu^s$; correspondingly, we
shall use the vector space of Radon measure on ${\mathbb R}^2$, which
we denote by ${\mathcal M}({\mathbb R}^2)$, for the reduced space of
vorticity functions, in place of the space of $H^{s-1}$ functions.

Recall that the co-Adjoint action of ${\mathcal D}_\mu^s$ on
$H^{s-1}({\mathbb R}^2)$ is given by
\begin{equation}
  \operatorname{Ad}^*_\eta(q) = q \circ \eta \,.
\end{equation}
We shall need to define the notion of weak co-adjoint action of
${\mathcal G}$ on ${\mathcal M}({\mathbb R}^2)$.  First, note that the
operation $\operatorname{Ad}^*\colon {\mathcal G} \times {\mathcal
M}({\mathbb R}^2) \rightarrow {\mathcal M}({\mathbb R}^2)$ given by
$\operatorname{Ad}^*_\eta (q) = q \circ\eta$ is well-defined.  Next,
define the {\it weak} co-Adjoint action $\operatorname{w-Ad}^* \colon
{\mathcal G} \times {\mathcal M}({\mathbb R}^2) \rightarrow {\mathcal
M}({\mathbb R}^2)$ by
\begin{equation}
  \int_{\mathbb R} \int_{{\mathbb R}^2}
    \operatorname{w-Ad}^*_\eta(q) \cdot \phi\, \d x \, \d t 
  = \int_{\mathbb R}\int_{{\mathbb R}^2} 
    q \cdot (\phi \circ \eta) \, \d x \, \d t
\end{equation}
for all $\phi \in C^\infty_0({\mathbb R}\times {\mathbb R}^2)$.

It follows that if $\eta_t$ is a $C^1$ curve in ${\mathcal G}$ such
that $e=\eta_0$ and $u= (d/dt)|_{t=0} \eta_t$, then we may---computing
the time derivative of $\operatorname{w-Ad}^*_{\eta_t}(q)$ at
$t=0$---define the weak analogue of the algebra co-adjoint action by
\begin{equation}
  \int_{\mathbb R} \int_{{\mathbb R}^2} 
    \operatorname{w-ad}^*_u ( q)\cdot \phi \, \d x \, \d t 
  = \int_{\mathbb R}\int_{{\mathbb R}^2} q \cdot
    \bigl( 
      \partial_t \phi + u\cdot\operatorname{grad}\phi
    \bigr) \, \d x\, \d t 
\end{equation}
for every $\phi\in C^\infty_0({\mathbb R} \times \R^2)$.  Recall that
the classical co-adjoint action is defined by
$\operatorname{ad}^*_{u_t} ( q_t)= (d/dt)|_{t=0}(\eta_t^* q_t)$ where
$(d/dt)|_{t=0}(\eta_t^* q_t) = \partial_t q_t + \mathcal{L}_{u_t}
q_t$.  In two dimensions, the Lie derivative term $\mathcal{L}_{u_t}
q_t$ reduces to $u_t \cdot \operatorname{grad} q_t$.

\begin{theorem}
For any $q_0 \in {\mathcal M}({\mathbb R}^2)$, let ${\mathcal
O}_{q_0}$ denote the co-Adjoint orbit $\{q \colon q= q_0 \circ \eta,
\eta\in {\mathcal G}\}$.  The weak co-adjoint action of
$C^0_{\mathrm{div}}({\mathbb R}^2)$ on ${\mathcal M}( {\mathbb R}^2)$
is well-defined, and solutions of the second-grade fluids equations or
of Chorin's vortex blob method with initial data $q_0$ leave
${\mathcal O}_{q_0}$ invariant.
\end{theorem}
\begin{proof}
The result immediately follows from the fact that the vanishing of the
weak co-adjoint action is equivalent to the weak formulation of
(\ref{e.alpha-vorticity}).  Theorem~\ref{t.well-posedness}, giving
global well-posedness of weak solutions, then concludes the argument.
\end{proof}

\section{Convergence}

We can now prove convergence of the flow of the vortex blob method to
the flow of the Euler equations.  This is done in two steps.  First we
show that the averaged Euler equation, or vortex method PDE,
approximates the Euler equation as $\alpha \to 0$ for bounded
vorticity fields.  In the second step, we prove that continuous
solutions of the averaged Euler equation can be approximated by
measure-valued ones.  These two results together imply
Theorem~\ref{t.convergence}.

\begin{lemma} \label{l.difference}
  Let $q_0 \equiv \omega_0 \in \l1(\R^2) \cap \l\infty(\R^2)$.  Then
  for every $T>0$ there exists a positive constant $C(T)$ such that
\begin{equation}
  \sup_{t \in [0,T]} \sup_{x \in \R^2}
    \bigl| \eta^\alpha (x,t) - \eta(x,t) \bigr|
  \leq C(T) \, \alpha^{\e^{-T}} \,.
\end{equation}
\end{lemma}

\begin{proof}
  We estimate the difference of the Euler and Euler-$\alpha$ flow
  maps:
\begin{align}
  \bigl| \eta^\alpha (x,t) - \eta (x,t) \bigr| 
  & \leq \int_0^t
         \bigl|
           u^\alpha \circ \eta^\alpha - u \circ \eta
         \bigr| \, \d s \notag \\
  & \leq \int_0^t \int_{\R^2} 
         \bigl| 
           K^\alpha (\eta^\alpha(x),\eta^\alpha(y)) -
           K        (\eta^\alpha(x),\eta^\alpha(y))
         \bigr| \, |\omega_0(y)| \,
         \d y \, \d s \notag \\
  & \quad +
         \int_0^t \int_{\R^2} 
         \bigl| 
           K (\eta^\alpha(x),\eta^\alpha(y)) -
           K (\eta(x),\eta^\alpha(y))
         \bigr| \, |\omega_0(y)| \,
         \d y \, \d s \notag \\
  & \quad +
         \int_0^t \int_{\R^2}
         \bigl| 
           K (\eta(x),\eta^\alpha(y)) -
           K (\eta(x),\eta(y))
         \bigr| \, |\omega_0(y)| \,
         \d y \, \d s \notag \\
  & \equiv
         \int_0^t (I_1 + I_2 + I_3) \, \d s
  \label{e.4} 
\end{align}
To estimate $I_1$, we note that on $\R^2$, the difference of the
kernels is explicitly given by $K^\alpha (r) - K (r) = K_1
(r/\alpha)/\alpha$, so that
\begin{align}
  I_1 & =    \int_{\R^2} 
               \frac1\alpha \, K_1 \Bigl( \frac{|x-y|}\alpha \Bigr)
               \, |\omega^\alpha (y,s)| \, \d y \notag \\
      & \leq 2 \pi
             \int_0^\infty
               \frac1\alpha \, K_1 \Bigl( \frac{r}\alpha \Bigr)
               r \, \d r \,
             \norm{\omega^\alpha (s)}{\l\infty} \notag \\
      & \leq c \, \alpha \, \norm{\omega_0}{\l\infty}
\end{align}
The other two integrals can be estimated by using the quasi-Lipschitz
conditions, Lemma~\ref{l.quasi-lipshitz.1} and
Lemma~\ref{l.quasi-lipshitz.2}, respectively.  One finds that
\begin{equation}
  I_2
  \leq c \,
       \varphi \bigl( \eta^\alpha(x,s) - \eta(x,s) \bigr) \,
       \bigl(
         \norm{\omega_0}{\l1} + \norm{\omega_0}{\l\infty}
       \bigr) \,,
\end{equation}
and
\begin{align}
  I_3
  & =    \int_{\R^2}
         \bigl|
           K (\eta (x,s),y) -
           K (\eta (x,s), \eta^\alpha \circ \eta^{-1} (y,s))
         \bigr| \,
         |\omega(y,s)| \, \d y \notag \\
  & \leq c \, \sup_{x \in \R^2} \varphi
         \bigl(
           x - \eta^\alpha \circ \eta^{-1} (x,s)
         \bigr) \, 
         \bigl(
           \norm{\omega(s)}{\l1} + \norm{\omega(s)}{\l\infty}
         \bigr) \notag \\
  & =    c \, \sup_{x \in \R^2}  
         \varphi \bigl( \eta(x,s) - \eta^\alpha(x,s) \bigr) \,
         \bigl(
           \norm{\omega_0}{\l1} + \norm{\omega_0}{\l\infty}
         \bigr) \,.
\end{align}
By inserting the bounds for $I_1$ to $I_3$ back into \eqref{e.4} and
taking the supremum on both sides, we obtain the log-Gronwall
inequality
\begin{equation}\label{logG}
  \sup_{x \in \R^2}
    \bigl| \eta^\alpha (x,t) - \eta (x,t) \bigr|
  \leq \int_0^t
       \Bigl[
         \alpha \, K_1 +
         K_2 \, \sup_{x \in \R^2}  
           \varphi \bigl( \eta^\alpha(x,s) - \eta(x,s) \bigr)
       \Bigr] \, \d s \,.
\end{equation}
To obtain explicit bounds that are valid on any finite interval of time
$[0,T]$, we set
\begin{equation}
  \rho (t) = \sup_{x \in \R^2}
    \bigl| \eta^\alpha (x,t) - \eta (x,t) \bigr| \,,
\end{equation}
and use the tangent approximation of the concave function $\varphi$;
namely, for any $\varepsilon \in (0,1)$,
\begin{equation}\label{concave}
  \varphi(r) \leq \varphi(\varepsilon) + \varphi'(\varepsilon) \, r 
             = (-\ln \varepsilon) \, r + \varepsilon \,.
\end{equation}
This makes the right-hand-side of (\ref{concave}) linear in $r$.  For
notational simplicity, we also rescale $\alpha$ and $t$ such that
$K_1=K_2=1$.  We substitute (\ref{concave}) into (\ref{logG}) and
obtain the usual Gronwall inequality; it follows that $\rho$ must
satisfy the differential inequality
\begin{equation}\label{diff.in}
  \dot\rho \leq (-\ln \varepsilon) \, \rho + \varepsilon + \alpha \,,
  \qquad
  \rho (0) = 0 \,.
\end{equation}
Setting $\varepsilon = \e^{-1} \alpha^{\exp(-t)}$ and integrating 
(\ref{diff.in}) with this choice of $\varepsilon (\alpha)$, we
find that
\begin{equation}
  \rho(t)
    \leq \frac{\e^t - 1}{\e} \, \alpha^{\e^{-t}}
         + \e^t \, \frac{\alpha^{\e^{-t}} - \alpha}{-\ln\alpha} \,. 
\end{equation}
Thus, $\rho = O(\alpha^{\exp(-T)})$ uniformly on $[0,T]$.
\end{proof}

In the following we will consider $\alpha$ as fixed and approximate
continuous data by measure valued data.  Let $\eta$, $q$, and $u$
denote quantities corresponding to a solution of the Euler-$\alpha$
equation with initial data $q_0 \in L^\infty(\R^2)$, and let $\eta^n$, $q^n$,
and $u^n$ denote a sequence of solutions to the Euler-$\alpha$
equation with initial data $q_0^n \in \mathcal{M} (\R^2)$ for every $n
\in \N$.  Then the following is true.

\begin{lemma} \label{l.difference.2}
Let $q_0 \in \l1(\R^2) \cap L^\infty(\R^2)$, and suppose that $q_0$ is
approximated by a sequence of measures in $\mathcal{M}(\R^2)$ such
that $q_0^n \rightharpoonup q_0$ weakly in $\mathcal{M}(\R^2)$, and
$\| q_0^n \|_{\mathcal M} \to \| q_0 \|_{\l1}$.  Then, for every
$T>0$,
\begin{equation}
  \lim_{n \to \infty}
  \sup_{t \in [0,T]} \sup_{x \in \R^2}
    \bigl| \eta^n (x,t) - \eta (x,t) \bigr| = 0 \,.
\end{equation}
\end{lemma}

\begin{proof}
As in the proof of Lemma~\ref{l.difference}, we estimate
\begin{align}
  \bigl| \eta & (x,t) - \eta^n (x,t) \bigr| \notag \\
  & \leq \int_0^t
         \biggl|
           \int_{\R^2}
             \bigl[
               K^\alpha (\eta(x,s), \eta(y,s)) \, q_0(y) -
               K^\alpha (\eta^n (x,s), \eta^n (y,s)) \, q_0^n (y)
             \bigr] \,
             \d y
         \biggr| \, \d s \notag \\
  & \leq \int_0^t \biggl| \int_{\R^2}
         K^\alpha (\eta(x), \eta(y)) \,         
         \bigl(
           q_0(y) - q_0^n(y)
         \bigr) \, \d y \biggr| \, \d s \notag \\
  & \quad + \int_0^t \int_{\R^2}
         \bigl|
           K^\alpha (\eta(x), \eta(y)) -
           K^\alpha (\eta(x), \eta^n(y))
         \bigr| \,
         \lvert
           q_0^n(y)
         \rvert \, \d y \, \d s \notag \\
  & \quad + \int_0^t \int_{\R^2}
         \bigl|
           K^\alpha (\eta(x), \eta^n(y)) -
           K^\alpha (\eta^n(x), \eta^n(y))
         \bigr| \,
         \lvert
           q_0^n(y)
         \rvert \, \d y \, \d s \notag \\
  & \equiv
         \int_0^t (J_1 + J_2 + J_3) \, \d s
  \label{e.j}
\end{align}
We find, after a change of variables, that
\begin{equation}
  \sup_{x \in \R^2} J_1
  = \sup_{x \in \R^2}
    \left|
      \int_{\R^2}
        K^\alpha (x,y) \, (q_0-q_0^n)(\eta^{-1}(y)) \, \d y
    \right| \,.
\end{equation}
By Lemma~\ref{l.uniformity} below with $\phi(x-y) = K^\alpha (x,y)$
and $q_n(y) = (q_0-q_0^n)(\eta^{-1}(y))$, this expression converges to
zero as $n \to \infty$.  Moreover, by Lemma~\ref{l.quasi-lipshitz.3},
\begin{align}
  I_2 & \leq \sup_{x \in \R^2} \sup_{y \in \R^2}
             \bigl|
               K^\alpha (\eta(x), \eta(y)) -
               K^\alpha (\eta(x), \eta^n(y))
             \bigr| \,
             \norm{q_0^n}{\mathcal M} \notag \\
      & \leq \sup_{y \in \R^2}
             \frac{c}\alpha \, \varphi
             \biggl(
               \frac{\eta(y) - \eta^n(y)}\alpha
             \biggr) \,
             \norm{q_0^n}{\mathcal M}
\intertext{and}
  I_3 & \leq \frac{c}\alpha \, \varphi
             \biggl(
               \frac{\eta(x) - \eta^n(x)}\alpha
             \biggr) \,
             \norm{q_0^n}{\mathcal M} \,.         
\end{align}
By inserting these estimates back into \eqref{e.j} and taking the
supremum in $x$ on both sides, we obtain an integral inequality that
can be solved with the log-Gronwall inequality exactly as in the proof
of Lemma~\ref{l.difference}.  The result then follows.
\end{proof}

\begin{lemma} \label{l.uniformity}
Let $\{q_n\}$ be a sequence of measures in $\mathcal{M} (\R^2)$
converging weakly to zero with uniformly bounded total variation and
uniform decay at infinity.  Further assume that $\phi$ is a continuous
test function with $\phi \to 0$ as $\lvert x \rvert \to \infty$.  Then
\begin{equation}
  \lim_{n \to \infty} \sup_{x \in \R^2}
  \int_{\R^2} \phi (x-y) \, q_n(y) \, \d y = 0 \,.  
\end{equation}
\end{lemma}

\begin{proof}
Set $M = \sup_n \lVert q_n \rVert_{\mathcal{M}}$ and $M' = \sup_x
\lvert \phi(x) \rvert$.  Let $\varepsilon>0$ be fixed.  By assumption
on the $\{q_n\}$, there exists an $R>0$ such that for every $n \in
\N$,
\begin{equation}
  \int_{|y|>R} \lvert q_n(y) \rvert \, \d y
  < \frac\varepsilon{6M'} \,.
\end{equation}
Moreover, there exists an $R'>0$ such that $\lvert \phi(x) \rvert <
\epsilon/(2M)$ for $\lvert x \rvert > R'$.  Since $\phi$ is uniformly
continuous on compact sets, there exists $\delta>0$ such that $\lvert
\phi(x) - \phi(x') \rvert < \varepsilon/(3M)$ for all $x, x' \in B (0,
2R+R')$ with $\lvert x - x' \rvert < \delta$.  Cover $B (0, R+R')$
with finitely many balls of radius $\delta$ and denote the centers of
these balls by $x_i$, $i \in I$.  Choose $N$ large enough such that
for $n \geq N$,
\begin{equation}
  \max_{i \in I}
  \biggl|
    \int_{\R^2} \phi(x_i-y) \, q_n(y) \, \d y
  \biggr|
  < \frac\varepsilon3 \,.
\end{equation}
Then for $\lvert x \rvert < R+R'$ there exists an $i \in I$ such that
$\lvert x - x_i \rvert < \delta$, and
\begin{align}
  \biggl|
    \int_{\R^2} & \phi(x-y) \, q_n(y) \, \d y
  \biggr| \notag \\
  & \leq \int_{\R^2}
            \bigl| \phi(x-y) - \phi(x_i-y) \bigr| \, |q_n(y)| \, \d y
         +
         \biggl|
           \int_{\R^2} \phi(x_i-y) \, q_n(y) \, \d y
         \biggr| \notag \\
  & \leq \sup_{\lvert y \rvert \leq R}
           \bigl| \phi(x-y) - \phi(x_i-y) \bigr| \, 
           \norm{q_n}{\mathcal M}
         + 2 \sup_{x \in \R^2} \lvert \phi(x) \rvert
           \int_{|y|>R} |q_n(y)| \, \d y
         + \frac\varepsilon3 \notag \\
  & \leq \frac\varepsilon{3M} \, M +
         2 M' \, \frac{\varepsilon}{6M'} +
         \frac\varepsilon3
    =    \varepsilon \,.         
\end{align}
On the other hand, if $\lvert x \rvert \geq R+R'$, then
\begin{align}
  \biggl|
    \int_{\R^2} & \phi(x-y) \, q_n(y) \, \d y
  \biggr| \notag \\
  & \leq
    \sup_{|x| \geq R'} \lvert \phi(x) \rvert
    \int_{|y|<R} |q_n(y)| \, \d y
  + \sup_{x \in \R^2} \lvert \phi(x) \rvert
    \int_{|y|<R} |q_n(y)| \, \d y \notag \\
  & \leq \frac\varepsilon{2M} \, M +
    M' \,\frac\varepsilon{6M'}
  < \varepsilon \,.
\end{align}
This completes the proof.
\end{proof}

\section*{Acknowledgments}
The authors thank Thomas Beale, David Levermore, Christian Lubich,
Jerry Marsden, Tudor Ratiu, and Edriss Titi for many interesting
discussions.  SS was partially supported by the NSF-KDI grant
ATM-98-73133.  MO was partially supported by the SFB 382 of the German
Science Foundation.

\end{document}